\documentclass[conference, letterpaper, final]{IEEEtran}
\IEEEoverridecommandlockouts

\usepackage{amsmath,amssymb,amsfonts}
\usepackage{cite}
\usepackage[linesnumbered, ruled, vlined]{algorithm2e}
\SetKwBlock{DoParallel}{do in parallel}{end}
\SetKwFor{ForParallel}{for}{do in parallel}{end for}
\usepackage{graphicx}
\usepackage{textcomp}
\usepackage{xcolor}
\usepackage{paralist} 
\usepackage[colorlinks,linkcolor=black,bookmarksopen,bookmarksnumbered,citecolor=black,urlcolor=black]{hyperref}
\usepackage{cleveref}
\Crefname{definition}{Definition}{Definitions}
\Crefname{proposition}{Proposition}{Propositions}
\Crefname{theorem}{Theorem}{Theorems}
\Crefname{figure}{Fig.}{Figs.}
\Crefname{equation}{Eq.}{Eqs.}
\Crefname{section}{Section}{Sections}
\Crefname{subsection}{Section}{Sections}
\Crefname{subsubsection}{Section}{Sections}
\Crefname{algorithm}{Algorithm}{Algorithms}
\usepackage{tikz}
\usetikzlibrary{positioning}
\usepackage{float}
\usepackage{booktabs}
\usepackage{multirow}
\usepackage{etoolbox}
\makeatletter
\patchcmd{\@makecaption}
  {\scshape}
  {}
  {}
  {}
\makeatother

\newcommand{\orcid}[1]{\href{https://orcid.org/#1}{\includegraphics[width=0.6em]{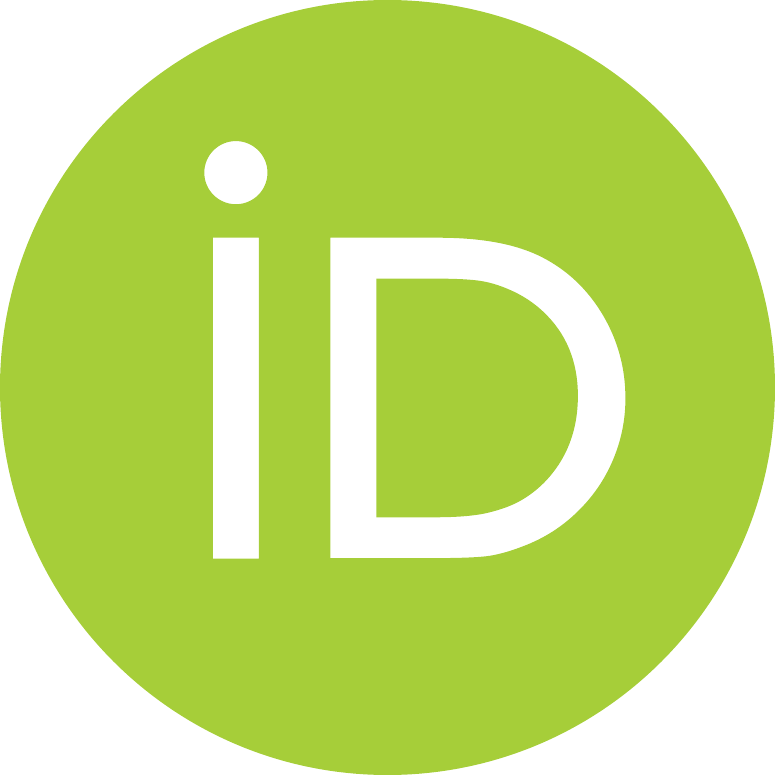}}}

\usepackage{glossaries}
\glsdisablehyper
\newacronym{oc}{OC}{optimal control}
\newacronym{mpc}{MPC}{model predictive control}
\newacronym{admm}{ADMM}{alternating direction method of multipliers}
\newacronym{pcg}{PCG}{preconditioned conjugate gradient}
\newacronym{lqr}{LQR}{linear-quadratic regulator}
\newacronym{kkt}{KKT}{Karush-Kuhn-Tucker}
\newacronym{ddp}{DDP}{differential dynamic programming}
\newacronym{nlp}{NLP}{nonlinear programming}
\newacronym{lu}{LU}{lower-upper}
\newacronym{fddp}{FDDP}{feasibility-driven DDP}
\newacronym{com}{CoM}{Center of Mass}

\usepackage{amssymb}
\usepackage{amsmath}
\usepackage{mathabx}
\usepackage{mathtools}

\newcommand{\smallddots}{\scalebox{0.7}{$\ddots$}}
\newcommand{\smallvdots}{\scalebox{0.7}{$\vdots$}}
\newcommand{\smallcdots}{\scalebox{0.7}{$\cdots$}}

\newcommand{\dynFunc}[1][]{\mathbf{f}_{#1}}
\newcommand{\consFunc}[1][]{\mathbf{h}_{#1}}
\newcommand{\termFunc}[1][]{\mathbf{r}_{#1}}

\newcommand{\grad}[1][]{\boldsymbol{\bigtriangledown_{#1}}}
\newcommand{\defeq}{\overset{\Delta}{=}}

\newcommand{\nx}{{n_{x}}} %
\newcommand{\ntau}{{n_{u}}} %
\newcommand{\nv}{{n_{v}}} %
\newcommand{\nc}{{n_{c}}} %
\newcommand{\nr}{{n_{r}}} %
\newcommand{\nw}{{n_{w}}} %
\newcommand{\nb}{{n_{b}}} %

\newcommand{\state}[1][]{\mathbf{x}_{#1}} %
\newcommand{\ctrl}[1][]{\mathbf{u}_{#1}} %
\newcommand{\termpoint}[1][]{\mathbf{c}_{#1}} %
\newcommand{\bp}[1][]{\mathbf{p}_{#1}} %
\newcommand{\stateSeq}{\state[s]}
\newcommand{\ctrlSeq}{\ctrl[s]}
\newcommand{\stateone}{\mathbf{x}_{1}}
\newcommand{\statetwo}{\mathbf{x}_{2}}
\newcommand{\stateN}{\state[N]}

\newcommand{\wrch}[1][]{\boldsymbol{\lambda}_c}

\newcommand{\initState}{\state[0]} %
\newcommand{\initCond}{\tilde{\state}_{0}}

\newcommand{\valFunc}[1][]{\mathcal{V}_{#1}}
\newcommand{\valFuncHat}[1][]{\hat{\mathcal{V}}_{#1}}
\newcommand{\valFuncCheck}[1][]{\widecheck{\mathcal{V}}_{#1}}
\newcommand{\qualFunc}[1][]{\mathbf{Q}_{#1}}
\newcommand{\costGrad}[1][]{\boldsymbol{\ell}_{#1}}
\newcommand{\lagHess}[1][]{\mathcal{L}_{#1}}
\newcommand{\mulpDyn}[1][]{\boldsymbol{\xi}_{#1}}
\newcommand{\mulpDynNext}[1][]{\boldsymbol{\xi}_{#1}^{\mathbf{+}}}

\newcommand{\mulpEqNext}[1][]{\boldsymbol{\gamma}_{#1}^{\mathbf{+}}}

\newcommand{\ffPolicy}[1][]{\mathbf{k}_{#1}}

\newcommand{\fbPolicy}[1][]{\mathbf{K}_{#1}}

\newcommand{\dynFeas}[1][]{\mathbf{\bar{f}}_{#1}}

\newcommand{\termState}{\state[N]} %
\newcommand{\cost}[1][]{\ell_{#1}} 

\newcommand{\dx}[1][]{\delta\state[#1]}
\newcommand{\du}[1][]{\delta\ctrl[#1]}

\newcommand{\dxHat}[1][]{\delta\hat{\state}_{#1}}
\newcommand{\duHat}[1][]{\delta\hat{\ctrl}_{#1}}
\newcommand{\dxNextHat}[1][]{\delta\hat{\mathbf{\state}}_{#1}{\mathbf{'}}}
\newcommand{\ffPolicyHat}[1][]{\hat{\boldsymbol{\pi}}_{#1}}
\newcommand{\fbPolicyHat}[1][]{{\boldsymbol{\Pi}}_{#1}}
\newcommand{\dX}[1][]{\delta\widecheck{\mathbf{X}}_{#1}}
\newcommand{\dU}[1][]{\delta\widecheck{\mathbf{U}}_{#1}}
\newcommand{\dL}[1][]{\delta\widecheck{\mathbf{\Lambda}}_{#1}}

\newcommand{\consFeas}[1][]{\mathbf{\bar{h}}_{#1}}
\newcommand{\termFeas}[1][]{\mathbf{\bar{r}}_{#1}}
\newcommand{\mulConsNext}{\boldsymbol{\beta}^{+}}

\newcommand{\mulpDynNextHat}[1][]{\boldsymbol{\widehat{\xi}}_{#1}^{\mathbf{+}}}

\newcommand{\mulpEqNextHat}[1][]{\boldsymbol{\widehat{\gamma}}_{#1}^{\mathbf{+}}}
\newcommand{\ffPolicyCheck}[1][]{\widecheck{\boldsymbol{\pi}}_{#1}}
\newcommand{\fbPolicyEq}[1][]{\mathbf{\Psi}_{#1}}
\newcommand{\qualFuncEq}[1][]{\tilde{\mathbf{Q}}_{#1}}
\newcommand{\transpose}{\intercal}
\newcommand{\zeroVec}[1][]{\mathbf{0}_{#1}}
\newcommand{\eyeMatrix}[1][]{\mathbf{I}_{#1}}

\newcommand{\R}[1][]{\mathbb{R}^{#1}}
\newcommand{\Spp}[1][]{S^{#1}_{++}}

\newcommand{\bA}[1][]{\mathbf{A}_{#1}}
\newcommand{\bB}[1][]{\mathbf{B}_{#1}}

\newcommand{\bWcheck}{{\widecheck{\mathbf{W}}}}

\newcommand{\bw}[1][]{\mathbf{w}_{#1}}
\newcommand{\bwHat}{\widehat{\mathbf{w}}}

\newcommand{\by}[1][]{\mathbf{y}_{#1}}

\newcommand{\ba}[1][]{\mathbf{a}_{#1}}
\newcommand{\bb}[1][]{\mathbf{b}_{#1}}

\newcommand{\schur}[1][]{\mathbf{S_{#1}}}
\newcommand{\nullspace}[1][]{\bar{\mathbf{Z}}_{#1}}
\newcommand{\rangespace}[1][]{\bar{\mathbf{Y}}_{#1}}
\newcommand{\nullspaceEq}[1][]{\mathbf{Z}_{#1}}
\newcommand{\rangespaceEq}[1][]{\mathbf{Y}_{#1}}
\newcommand{\nullVector}[1][]{\bar{\mathbf{z}_{#1}}}
\newcommand{\rangeVector}[1][]{\bar{\mathbf{y}_{#1}}}
\newcommand{\nullEqVector}[1][]{\mathbf{z}_{#1}}
\newcommand{\rangeEqVector}[1][]{\mathbf{y}_{#1}}
\newcommand{\nz}[1][]{{n_{\bar{Z}_{#1}}}}
\newcommand{\ny}[1][]{n_{\bar{Y}_{#1}}}
\newcommand{\rs}[1][]{\bar{y}}
\newcommand{\ns}[1][]{\bar{z}}

\newcommand{\rev}[1]{\textcolor{black}{#1}}

\title{Endpoint-Explicit Differential Dynamic Programming via Exact Resolution}
\author{
    Maria Parilli\orcid{0009-0009-3210-3155}\quad
    Sergi Martinez\orcid{0009-0001-0027-8358}\quad
    Carlos Mastalli\orcid{0000-0002-0725-4279}
}

\begin{document}
\maketitle

\begin{abstract}
We introduce a novel method for handling endpoint constraints in constrained~\gls{ddp}.
Unlike existing approaches, our method guarantees quadratic convergence and is exact, effectively managing rank deficiencies in both endpoint and stagewise equality constraints.
It is applicable to both forward and inverse dynamics formulations, making it particularly well-suited for~\gls{mpc} applications and for accelerating~\gls{oc} solvers.
We demonstrate the efficacy of our approach across a broad range of robotics problems and provide a user-friendly open-source implementation within \textsc{Crocoddyl}.
\end{abstract}
\glsresetall

\section{Introduction}
Efficient and exact handling of endpoint constraints is a critical component for the parallelization of algorithms in~\gls{oc}.
In recent years, various approaches have emerged to address this challenge through endpoint constraints.
One straightforward approach is leveraging algorithms for distributed optimization such as the~\gls{admm}~\cite{boyd-admmbook}.
This approach was recently applied by Bishop et al.~\cite{bishop2024-reluqp}.
However,~\gls{admm} is inherently inexact and achieves only linear convergence at best, due to its reliance on augmented Lagrangian techniques \rev{(see \cite[Section II-C, p. 6]{yang2022-surveyadmm})}.%

A promising alternative for parallelizing computations is the use endpoint-explicit Riccati recursions, as proposed by Laine and Tomlin~\cite{laine2019-parallel_lqr}.
These recursions offer the potential for quadratic convergence.
However, they manage endpoint constraints through pseudo-inverses~\cite{laine2019-endpointlqr}, a process that is both computationally expensive and numerically unstable.
Consequently, this approach is limited to~\gls{lqr} frameworks and struggles with more complex, nonlinear~\gls{oc} problems.

In many real-time robotics applications, strict terminal \rev{constraints} must be met within tight computational time limits.
In this context,~\gls{ddp} has become the preferred approach, as it leverages the Markovian structure of the problem through Riccati recursions.
This results in significantly higher computational efficiency compared to traditional sparse linear solvers (e.g., MA27, MA57, MA97~\cite{HSL}) commonly used in~\gls{nlp}.
However, recent advances in constrained DDP methods (e.g.,~\cite{howell2019-altro,aoyama2021-alddp,pavlov2021-ipddp,jallet2022-constrainedddp}) have focused mainly on stagewise constraints, often treating endpoint constraints inexactly or using penalty-based methods, which can compromise performance in tasks requiring precise endpoint \rev{constraints}.
\begin{figure}[t]
    \centering
    \href{https://youtu.be/RBohdOhgbWw}{
    \begin{minipage}{0.48\linewidth}
        \centering
        \includegraphics[width=0.95\linewidth]{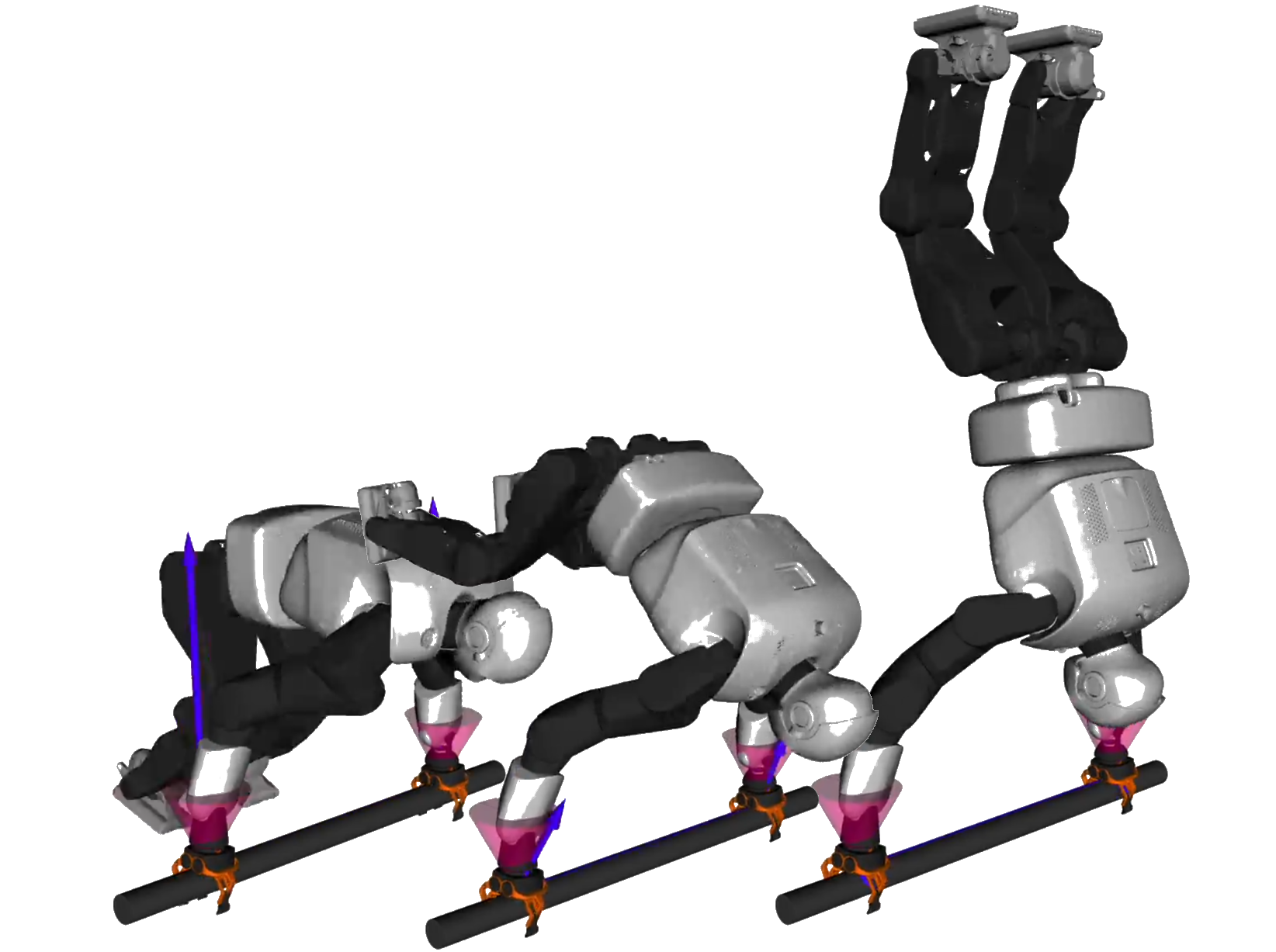} \\
        \includegraphics[width=0.95\linewidth]{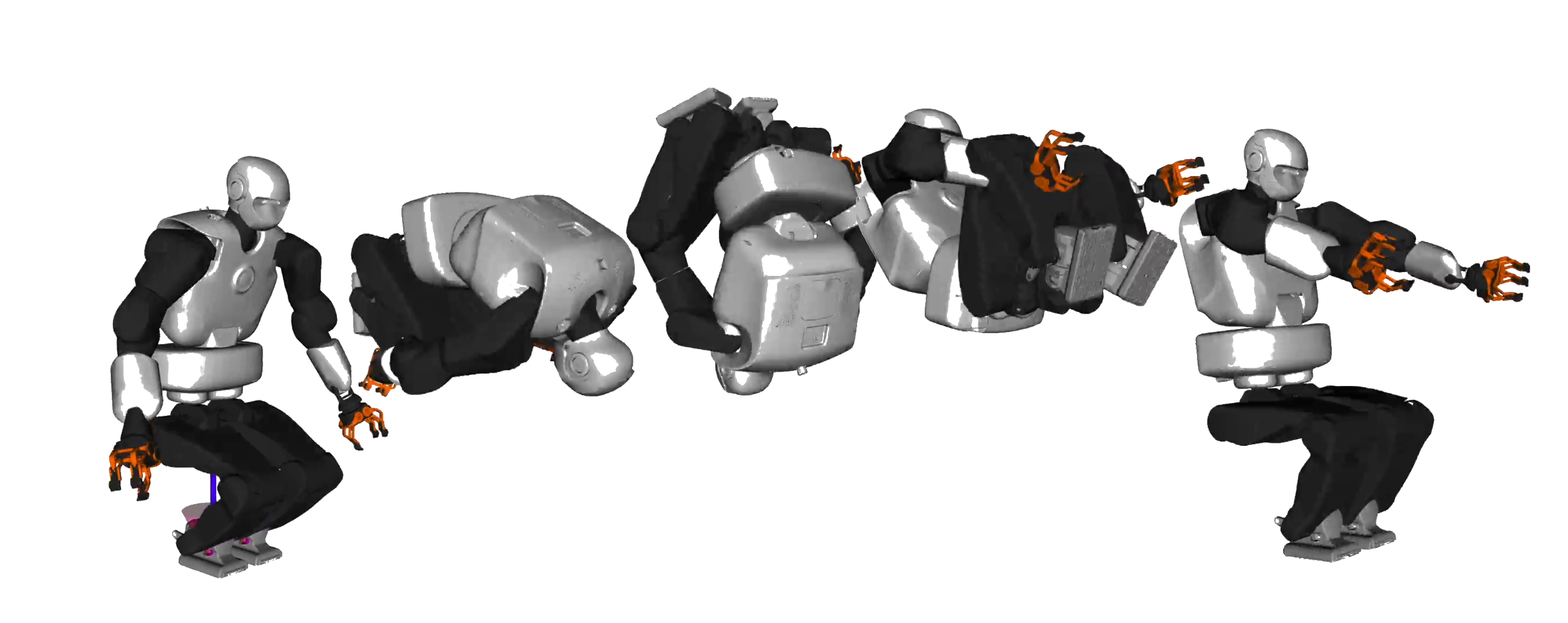}
    \end{minipage}
    \begin{minipage}{0.48\linewidth}
        \centering
        \includegraphics[width=0.95\linewidth]{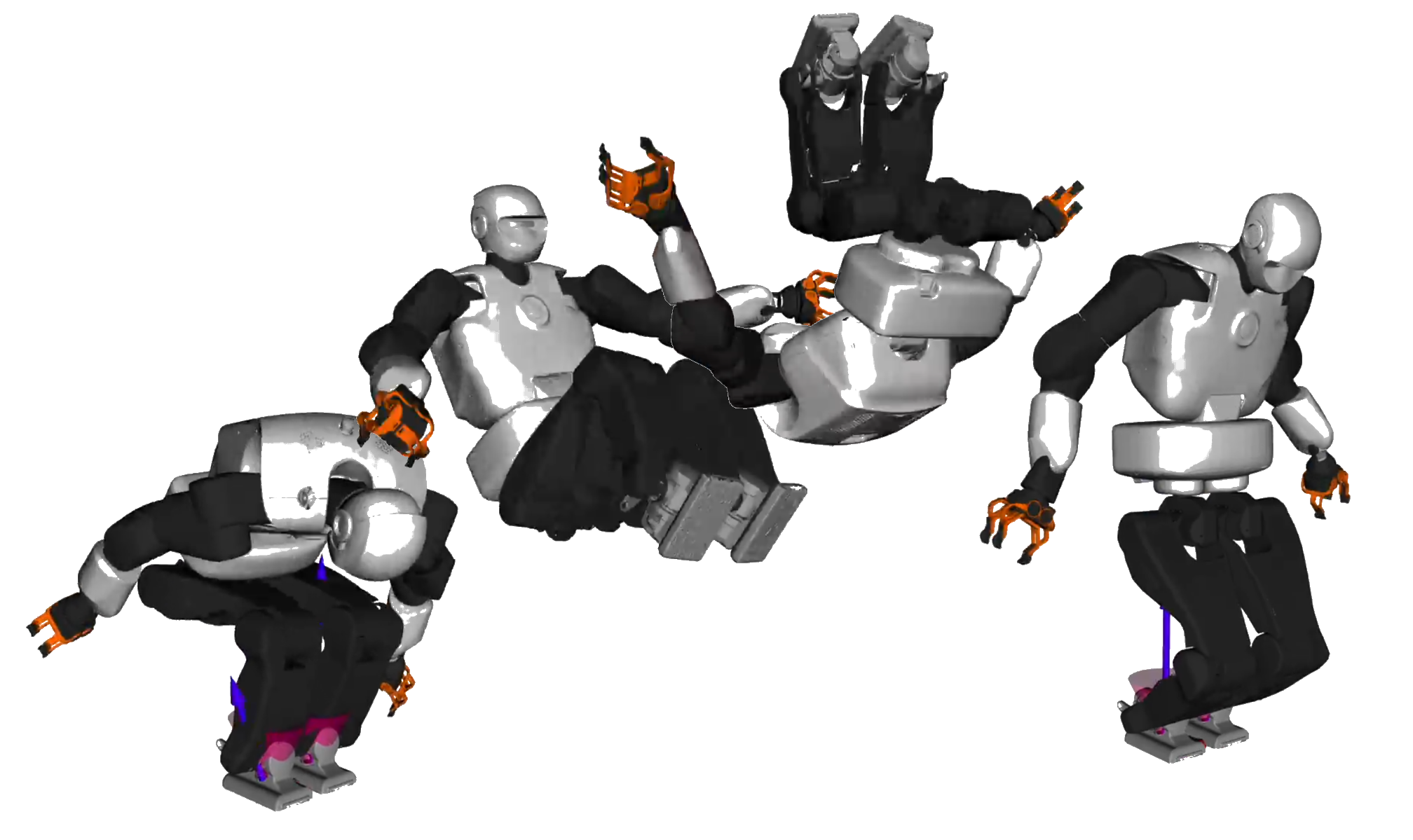} \\
        \includegraphics[width=0.95\linewidth]{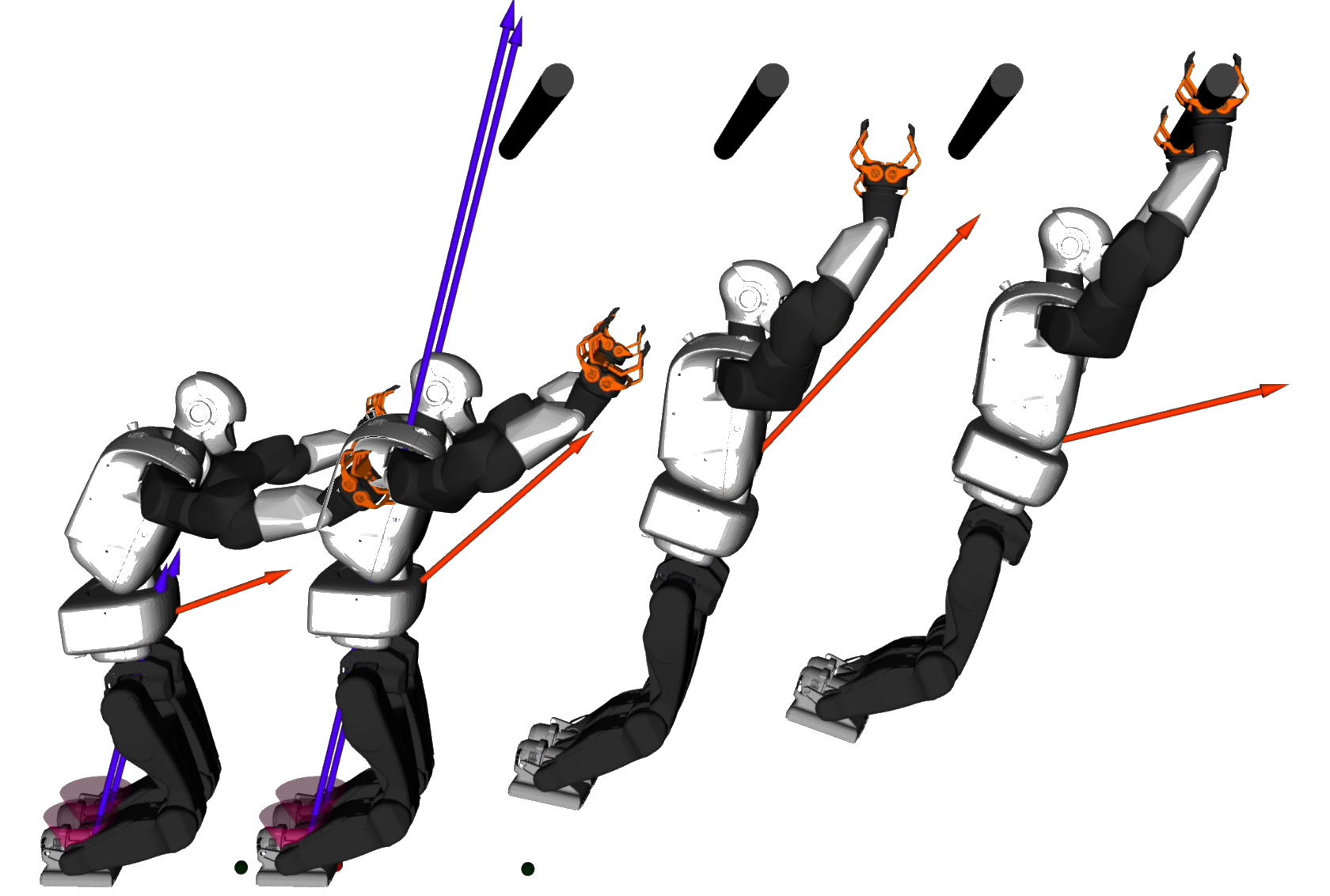}
    \end{minipage}
    }\vspace{-0.25em}
    \caption{Snapshots of Talos performing various gymnastic maneuvers computed using our endpoint-explicit~\gls{ddp} algorithm.
    The images show Talos executing a (top-left) handstand, (top-right) backflip, (bottom-left) frontflip, and (bottom-right) monkey bar maneuver.
    For each maneuver, the feet or hand placements were specified as endpoint constraints.
    \rev{To watch the video, click the picture or see \texttt{\url{https://youtu.be/RBohdOhgbWw}}}.}\vspace{-1em}
    \label{fig:snapshots}
\end{figure}

\subsection{Contribution}
In this work, we address the limitation of existing methods for handling endpoint constraints in~\gls{oc}.
Our approach is both exact and computationally efficient, offering quadratic convergence while managing rank-deficient endpoint constraints.
Furthermore, it accommodates nonlinearities in both forward and inverse dynamics formulations, making it highly suitable for~\gls{mpc} applications.

We demonstrate the effectiveness of our method across various robotics applications~(\Cref{fig:snapshots}). %
Our approach is adaptable to various stagewise constraints and readily accessible to practitioners through \textsc{Crocoddyl}~\cite{mastalli2020-crocoddyl}.

\section{Quadratic programs with endpoint constraints}\label{sec:abstract-solution}
Before diving into our algorithm, we first introduce a key abstraction that facilitates its development.
Specifically,~\gls{oc} problems subject to stagewise and endpoint constraints can be locally resolved by solving the following quadratic program:
\begin{equation}\label{eq:terminal-constrained-oc-abstraction}
\begin{aligned}
&\min_{\bw} \; \frac{1}{2} \bw^{\transpose} \bA \bw - \bw^{\transpose} \ba\\
& \text{subject to} \quad \bB \bw = \bb,
\end{aligned}
\end{equation}
at each Newton step to find the nonlinear roots that define the \gls{kkt} point~\cite{nocedal-optbook}.

In this abstraction, $\bw\in\R[\nw]$ represents the primal and dual decision variables encountered in an~\gls{oc} problem without endpoint constraints (e.g., search direction for states, controls, and multipliers of the dynamics constraints).
The matrix $\bA\in\Spp[\nw]$ is a symmetric, positive definite, banded, and large, while $\bB\in\R[\nb\times\nw]$ is the wide ($\nb\ll\nw$) Jacobian of the endpoint constraints, with $\bb\in\R[\nb]$ as its bias term.

The Lagrangian of~\Cref{eq:terminal-constrained-oc-abstraction} is defined as:
\begin{equation}
\label{eq:lagrangian-abstraction}
    \lagHess(\bw,\by) = \frac{1}{2}\bw^{\transpose}\bA\bw - \bw^{\transpose}\ba + \by^{\transpose}(\bB\bw -\bb),
\end{equation}
where $\by\in\R[\rev{n_b}]$ represents the Lagrange multipliers for the endpoint constraints.
This allows us to establish the \gls{kkt} conditions for~\Cref{eq:terminal-constrained-oc-abstraction} as follows
\begin{equation}
\label{eq:kkt-conditions-abstraction}
\begin{aligned}
   \grad[\bw]\lagHess(\bw,\by) = \bA\bw - \ba + \bB^{\transpose}\by = \zeroVec, & \quad \text{(primal feasibility)} \\
   \rev{\grad[\by]}\lagHess(\bw,\by) = \bB\bw - \bb = \zeroVec, & \quad\text{(dual feasibility)}
\end{aligned}
\end{equation}
which leads to the large symmetric saddle point system:
\begin{equation}\label{eq:saddle-point-system}
\begin{bmatrix}
\bA & \bB^{\transpose} \\
\bB &
\end{bmatrix} 
\begin{bmatrix}
\bw \\ \by 
\end{bmatrix}
= \begin{bmatrix}
\ba \\ \bb
\end{bmatrix}\in\R[\nw + \nb].
\end{equation}

Numerous methodologies have been proposed in the literature for factorizing such saddle point systems.
Here, we briefly overview two prominent techniques: the classical Schur-complement approach and nullspace factorization~\cite{rees2018-nullspace_saddle_point_systems}.
These methods are crucial due to their exactness, which is essential for ensuring fast convergence when for handling endpoint constraints.
This stands in contrast to inexact approaches, such as Krylov subspace method or relaxed proximal formulations~\cite{parish2014-proxalgorithms}.
The \rev{latter} can be viewed as a generalized Augmented Lagrangian method~\cite{gill2012-pdauglag}, where convergence depends on the accuracy of the estimated multiplier $\by[e]$.
Moreover, solving~\Cref{eq:saddle-point-system} leverages second-order derivatives to achieve quadratic convergence.

\subsection{Schur-complement resolution}\label{sec:abstract_schur_complement}
Since $\bA$ is square and nonsingular, we can solve \Cref{eq:saddle-point-system} as\vspace{-1.25em}
\begin{subequations}\label{eq:abstract_solution}
\begin{align}\label{eq:saddle_multipliers}
 &\by = -\schur^{-1}(\bb-\bB\bA^{-1}\ba),\\\label{eq:saddle_decision_var}
 &\bw = \bA^{ -1}\ba -\bA^{-1}\bB^{\transpose}\by,
\end{align}
\end{subequations}
where $\schur=\bB\bA^{-1}\bB^{\transpose}$ denotes the Schur complement, and its inversion can be efficiently computed via a sparse Cholesky decomposition.

From \Cref{eq:saddle_decision_var}, we observe that the decision vector $\bw$ decomposes into two terms:
\begin{equation}
\bw = \bwHat - \bWcheck \by,    
\end{equation}
where $\bwHat=\bA^{-1}\ba$ is the endpoint-independent solution, and $\bWcheck=\bA^{-1}\bB^{\transpose}\in\R[\nw\times\nb]$ accounts for the effect of the endpoint multiplier $\by$, i.e., the endpoint-dependent solution.

Finally, the banded structure of \(\bA\), which arises in optimal control due to its Markovian nature, allows for the application of Riccati recursions when solving \Cref{eq:abstract_solution}. This results in a highly efficient algorithm for solving arbitrary constraints.

\subsection{Nullspace resolution}\label{sec:abstract_nullspace}
If $\bB$ is rank deficient, the Schur complement of~\rev{\Cref{eq:abstract_solution}} becomes non-invertible, rendering the problem ill-posed.
\rev{A} straightforward solution is to parameterize the decision vector $\bw$ into its nullspace and range.
However, this strategy:
\begin{inparaenum}[(i\upshape)]
\item requires a costly nullspace decomposition due to the large and sparse nature of $\bB$, and
\item does not exploit the structure of $\bA$, preventing efficient computation of $\bwHat$ and $\bWcheck\by$.
\end{inparaenum}
Instead, we propose performing a nullspace decomposition on the Lagrange multiplier \(\by\).
Specifically, we parameterize \(\by\) as:
\begin{equation}\label{eq:decision_var_nullspace_parametrization}
\by=\nullspace\by[\nullVector]+\rangespace\by[\rangeVector],
\end{equation}
where $\nullspace\in\R[\rev{n_b}\times\nz]$ represents the nullspace basis of $\schur$ and $\rangespace\in\R[\rev{n_b}\times\ny]$ denotes its orthonormal basis.
This reduces~\Cref{eq:saddle-point-system} to the following saddle point system:
\begin{equation}\label{eq:reduced-saddle-point-system}
\begin{bmatrix}
\bB[\nullVector]\bA^{-1}\bB[\nullVector]^{\transpose} & \bB[\nullVector]\bA^{-1}\bB[\rangeVector]^{\transpose}\\ 
\bB[\rangeVector]\bA^{-1}\bB[\nullVector]^{\transpose} & \bB[\rangeVector]\bA^{-1}\bB[\rangeVector]^{\transpose}
\end{bmatrix}
\begin{bmatrix}
\by[\nullVector] \\ \by[\rangeVector] 
\end{bmatrix} = -
\begin{bmatrix}
\bb[\nullVector]-\bB[\nullVector]\bwHat \\ 
\bb[\rangeVector]-\bB[\rangeVector]\bwHat
\end{bmatrix},
\end{equation}
where $\bB[\nullVector]\defeq\nullspace^{\transpose}\bB$, $\bB[\rangeVector]\defeq\rangespace^{\transpose}\bB$, $\bb[\nullVector]\defeq\nullspace^{\transpose}\bb$ and $\bb[\rangeVector]\defeq\rangespace^{\transpose}\bb$.
Additionally, since $\bB[\nullVector]=\zeroVec$, the Lagrange multiplier is computed as:
\begin{equation}
\label{eq:saddle_decision_var_in_rangespace}
 \by = \rev{\rangespace\by[\rangeVector] =} -\rev{\rangespace}(\bB[\rangeVector]\bA ^{-1}\bB[\rangeVector]^{\transpose})^{-1}(\bb[\rangeVector] - \bB[\rangeVector]\bwHat).
\end{equation}
This formulation effectively handles rank \rev{deficiency} in $\bB$, allowing for efficient resolution of the system even when the Schur complement is \rev{invertible}.

\section{Optimal control with endpoint constraints}
Nonlinear~\gls{oc} problems with endpoint constraints can be formulated using either forward or inverse dynamics.
In inverse dynamics formulation, the problem \rev{involves} handling stagewise equality constraints as follows
\begin{equation}\label{eq:constrained-oc}
\begin{alignedat}{2}
&\min_{\stateSeq,\ctrlSeq} \; \cost[N](\termState) + \sum_{k=0}^{N-1} \cost[k](\state[k],\ctrl[k])\\
\text{s.t.} \quad 
&\initCond \ominus \initState = \zeroVec, \quad &\text{(initial \rev{constraint})} \\
&\dynFunc[k]({\state}_k,{\ctrl}_k) \ominus \state[k+1] = \zeroVec, \quad &\hspace{-3em}\text{(integrator / forward dyn.)}\\
&\consFunc[k](\state[k],\ctrl[k]) = \zeroVec, \quad &\text{(inverse dyn.)}\\
&\termFunc(\state[N]) = \zeroVec. \quad &\text{(endpoint \rev{constraint})}
\end{alignedat}
\end{equation}
The state $\state\in{\cal{X}}\subseteq\R[\nx]$ lies in a smooth manifold, with $\initCond$ specifying the system’s initial state. The control inputs $\ctrl[k]\in\R[\ntau]$ is optimized to minimize the objective, where:
\begin{inparaenum}[(i\upshape)]
\item $\cost[k]:\R[\nx]\times\R[\ntau]\rightarrow\R_{++}$ represents the stagewise cost for each time step,
\item $\cost[N]:\R[\nx]\rightarrow\R_{++}$ is the terminal cost,
\item $\dynFunc[k]:\R[\nx]\times\R[\ntau]\rightarrow\R[\nx]$ captures the system’s forward dynamics or integrator~\rev{as in~\cite{mastalli-inv_dyn_mpc}},
\item $\consFunc[k]:\R[\nx]\times\R[\ntau]\rightarrow\R[\nv+\nc]$ incorporates inverse dynamics and contact constraints, with $\nv$ representing the dimension of inverse dynamics and $\nc$ the dimension of contact constraints,
\item $\termFunc:\R[\nx]\rightarrow\R[\nr]$ defines the endpoint \rev{constraint}.
\end{inparaenum}
Moreover, $\ctrl$ represents joint efforts in forward dynamics formulations~\cite{mastalli2022-agilempc} or generalized accelerations and contact forces in inverse dynamics formulations~\cite{mastalli-inv_dyn_mpc}.

To solve this large optimization problem, we find the \gls{kkt} point using an iterative routine based on the Newton method.
This approach requires factorizing the following system of equations:
\begin{equation}
\label{eq:terminal-constrained-KKT}
\setlength{\arraycolsep}{0.7pt}
\begin{bmatrix}
& -\eyeMatrix \\
-\eyeMatrix & \smallddots\\
&& \lagHess[\state{\state[k]}] & \lagHess[\state{\ctrl[k]}] & \consFunc[{\state[k]}]^\transpose & \dynFunc[{\state[k]}]^\transpose\\
&& \lagHess[\ctrl{\state[k]}] & \lagHess[\ctrl{\ctrl[k]}] & \consFunc[{\ctrl[k]}]^\transpose & \dynFunc[{\ctrl[k]}]^{\transpose}\\
&& \consFunc[{\state[k]}] & \consFunc[{\ctrl[k]}]\\
&& \dynFunc[{\state[k]}] & \dynFunc[{\ctrl[k]}] && \smallddots\\
&&&&&& \costGrad[\state{\state[N]}] & \termFunc[{\state[N]}]^\transpose\\ 
&&&&&& \termFunc[{\state[N]}] & 
\end{bmatrix} 
\begin{bmatrix}
\mulpDynNext[0] \\ \smallvdots \\ \delta\state[k] \\ \delta\ctrl[k] \\ \mulpEqNext[k] \\\mulpDynNext[k+1] \\ \smallvdots \\ \dx[N] \\ \mulConsNext
\end{bmatrix}
= -
\begin{bmatrix}
\dynFeas[0] \\ \smallvdots \\ \costGrad[{\state[k]}] \\ \costGrad[{\ctrl[k]}] \\ \consFeas[k] \\ \dynFeas[k] \\ \smallvdots \\  \costGrad[{\state[N]}] \\ \termFeas
\end{bmatrix}
\end{equation}
In~\Cref{eq:terminal-constrained-KKT}, the infeasibilities are defined as follows: $\dynFeas[0]\coloneqq\initCond\ominus\initState$, $\consFeas[k]\coloneqq\consFunc[k](\state[k],\ctrl[k])$, $\dynFeas[k+1]\coloneqq\dynFunc[k](\state[k],\ctrl[k])\ominus\state[k+1]$, $\termFeas\coloneqq\termFunc(\state[N])$, which corresponds to initial \rev{constraint}, integrator, inverse dynamics, \rev{and endpoint constraint,} respectively.
Moreover, $\mulpDynNext[0]$, $\mulpEqNext[k]$, $\mulpDynNext[k]$, and $\mulConsNext$ are their Lagrange multipliers, where the $+$ notation stands for the next value, e.g., $\mulpDynNext\coloneqq\mulpDyn+\delta\mulpDyn$\rev{.}
Additionally, $\dynFunc[\bp]$, $\consFunc[\bp]$ represent the Jacobians of the dynamics and constraints, $\costGrad[\bp]$ denotes the cost gradient, where $\bp$ refers to either $\state$ or $\ctrl$.
The Lagrangian’s Hessian, $\lagHess[\bp\bp]$, is introduced below.

\subsection{Exploiting the structure of the problem}
\Cref{eq:terminal-constrained-KKT} matches the form introduced in~\Cref{eq:kkt-conditions-abstraction}.
Specifically, we can partition~\Cref{eq:terminal-constrained-KKT} as follows
\begin{equation}
\label{eq:constrained-saddle-matrices}
\begin{aligned}
&\bA = 
\setlength{\arraycolsep}{1.3pt}
\begin{bmatrix}
& -\eyeMatrix \\
-\eyeMatrix & \smallddots\\
&& \lagHess[\state{\state[k]}] & \lagHess[\state{\ctrl[k]}] & \consFunc[{\state[k]}]^\transpose & \dynFunc[{\state[k]}]^\transpose\\
&& \lagHess[\ctrl{\state[k]}] & \lagHess[\ctrl{\ctrl[k]}] & \consFunc[{\ctrl[k]}]^\transpose & \dynFunc[{\ctrl[k]}]^{\transpose}\\
&& \consFunc[{\state[k]}] & \consFunc[{\ctrl[k]}]\\
&& \dynFunc[{\state[k]}] & \dynFunc[{\ctrl[k]}] && \smallddots\\
&&&&&& \costGrad[\state{\state[N]}]
\end{bmatrix},\,\,
\ba = \rev{-}
\begin{bmatrix}
\dynFeas[0] \\ \smallvdots \\ \costGrad[{\state[k]}] \\ \costGrad[{\ctrl[k]}] \\ \consFeas[k] \\ \dynFeas[k] \\ \smallvdots \\  \costGrad[{\state[N]}]
\end{bmatrix},\\
&\bB = -
\setlength{\arraycolsep}{5.5pt}
\begin{bmatrix}
\zeroVec & \smallcdots & \zeroVec & \zeroVec & \zeroVec & \smallcdots & -\termFunc[{\state[N]}]
\end{bmatrix}, \quad
\bb = - \termFeas,
\end{aligned}
\end{equation}
where $\bA$ represents the~\gls{kkt} matrix associated to an optimal control problem without endpoint constraints, and $\ba$ denotes its~\gls{kkt} vector.

In~\Cref{sec:abstract-solution}, we initially outlined an abstract solution for this large saddle point system.
Now, we continue by detailing how to compute~\Cref{eq:abstract_solution}, leveraging the temporal structure of the problem for computational efficiency.

\subsection{Endpoint-independent search direction}\label{sec:unconstrained_search_dir}
To compute the term $\bwHat=\bA^{-1}\ba$, we can employ two factorization techniques: the Schur-complement and the nullspace approach.
Both methods enable the design of efficient Riccati recursions.
These Riccati recursions are carried out by factorizing the following system of equations backward in time:
\begin{equation}\label{eq:Riccati}
\rev{\begin{bmatrix}
\lagHess[\state\state] &\lagHess[\state\ctrl] & \consFunc[\state]^\transpose & \dynFunc[\state]^{\transpose}\\
\lagHess[\ctrl\state] & \lagHess[\ctrl\ctrl] & \consFunc[\ctrl]^\transpose & \dynFunc[\ctrl]^{\transpose}\\
\consFunc[\state] & \consFunc[\ctrl]\\
\dynFunc[\state] & \dynFunc[\ctrl] &&& -\eyeMatrix\\
&&& -\eyeMatrix & \valFunc[\state\state]'\\ 
\end{bmatrix}} 
\begin{bmatrix}
\dxHat\\ \duHat\\ \mulpEqNextHat\\ \mulpDynNextHat\\ \dxNextHat
\end{bmatrix}
= -
\begin{bmatrix}
\costGrad[\state]\\ \costGrad[\ctrl]\\ \consFeas\\ \dynFeas\\ \valFunc[\state]'
\end{bmatrix},
\end{equation}
where $\valFunc[\state{\state[N]}]$ and $\valFunc[{\state[N]}]$ correspond to $\costGrad[\state{\state[N]}]$ and $\costGrad[{\state[N]}]$, respectively. 
The notation $(\cdot)'$ refers to quantities at the next time step $k+1$, and $\hat{(\cdot)}$ denotes the endpoint-independent variables.
In each Riccati sweep, we first condense~\Cref{eq:Riccati} into the following system of equations
\begin{equation}\label{eq:qual-system}
\begin{aligned}
&\begin{bmatrix}
\qualFunc[\state\state] & \qualFunc[\state\ctrl] & \consFunc[\state]^\transpose\\
\qualFunc[\ctrl\state] & \qualFunc[\ctrl\ctrl] & \consFunc[\ctrl]^\transpose\\
\consFunc[\state] & \consFunc[\ctrl]
\end{bmatrix}
\begin{bmatrix}
\dxHat\\ \duHat\\ \mulpEqNextHat
\end{bmatrix}
=- 
\begin{bmatrix}
\qualFunc[\state]\\ \qualFunc[\ctrl]\\ \consFeas
\end{bmatrix}, \;
\end{aligned}
\end{equation}
where the $\qualFunc$'s terms are defined as follows
\begin{equation}
\label{eq:qual-funcs-def}
\begin{matrix}
\qualFunc[\state\state]\coloneqq\lagHess[\state\state]+\dynFunc[\state]^{\transpose}\valFunc[\state\state]'\dynFunc[\state], &
\qualFunc[\state]\coloneqq\costGrad[\state]+\dynFunc[\state]^{\transpose}\valFuncHat[\state]^+,\\
\qualFunc[\state\ctrl]\coloneqq\lagHess[\state\ctrl]+\dynFunc[\state]^{\transpose}\valFunc[\state\state]'\dynFunc[\ctrl], &
\qualFunc[\ctrl]\coloneqq\costGrad[\ctrl]+\dynFunc[\ctrl]^{\transpose}\valFuncHat[\state]^+,\\
\rev{\qualFunc[\ctrl\ctrl]}\coloneqq\lagHess[\ctrl\ctrl]+\dynFunc[\ctrl]^{\transpose} \valFunc[\state\state]'\dynFunc[\ctrl],
\end{matrix}
\end{equation}
with $\lagHess[\state\state]\coloneqq\costGrad[\state\state]+\valFuncHat[\state]'\cdot\dynFunc[\state\state]$, $\lagHess[\state\ctrl]\coloneqq\costGrad[\state\ctrl]+\valFuncHat[\state]'\cdot\dynFunc[\state\ctrl]$, $\lagHess[\ctrl\ctrl]\coloneqq\rev{\costGrad[\ctrl\ctrl]}+\valFuncHat[\state]'\cdot\dynFunc[\ctrl\ctrl]$, and $\valFuncHat[\state]^+\coloneqq\valFuncHat[\state]'+\valFunc[\state\state]'\dynFeas$.
Next, we compute $\duHat$ as a function of $\dxHat$, transforming~\Cref{eq:qual-system} into
\begin{equation}\label{eq:unconstrained-policy}
\begin{bmatrix}
\qualFunc[\ctrl\ctrl] & \consFunc[\ctrl]^\transpose \\
\consFunc[\ctrl]
\end{bmatrix}
\begin{bmatrix}
\duHat\\ \mulpEqNextHat
\end{bmatrix}
= - 
\begin{bmatrix}
\qualFunc[\ctrl] + \qualFunc[\ctrl\state]\dxHat \\
\consFeas + \consFunc[\state]\dxHat
\end{bmatrix}.
\end{equation}
By solving~\Cref{eq:unconstrained-policy} we obtain the endpoint-independent local control policy:
\begin{equation}\label{eq:control-policy}
\duHat= -\ffPolicyHat - \fbPolicyHat\dxHat,
\end{equation}
where $\ffPolicyHat$, $\fbPolicyHat$ are the feed-forward and feedback terms in~\gls{oc} problems without endpoint constraints.
Finally, we update its local approximation of the value function needed for the next Riccati sweep as follows
\begin{equation}\label{eq:value-func-approx-unconstrained}
\delta\valFuncHat=\frac{\Delta\valFuncHat[1]+2\Delta\valFuncHat[2]}{2}+\dxHat^{\transpose}(\valFuncHat[\stateone]+\valFuncHat[\statetwo])+\frac{1} {2}\dxHat^{\transpose}\valFunc[\state\state]\dxHat
\end{equation}
with terms computed as follows:
\begin{align}\label{eq:value-func-cont-unconstrained}\nonumber
&\Delta\valFuncHat[1]=\ffPolicyHat ^{\transpose}\qualFunc[\ctrl\ctrl]\ffPolicyHat,
&\Delta\valFuncHat[2]=-\ffPolicyHat^{\transpose}\qualFunc[\ctrl],\\\nonumber
&\valFuncHat[\stateone]=\fbPolicyHat^{\transpose}\qualFunc[\ctrl\ctrl] \ffPolicyHat-\qualFunc[\state\ctrl]\ffPolicyHat, &\valFuncHat[\statetwo]=\qualFunc[\state]-\fbPolicyHat^{\transpose}\qualFunc[\ctrl],\\
&\valFunc[\state\state]=\qualFunc[\state\state]-2\qualFunc[\state\ctrl]\fbPolicyHat+\fbPolicyHat^{\transpose}\qualFunc[\ctrl\ctrl]\fbPolicyHat.
\end{align}
For forward dynamics, \Cref{eq:value-func-cont-unconstrained} simplifies, where $\Delta\valFuncHat[1]=\ffPolicyHat\qualFunc[\ctrl]$, $\valFuncHat[\stateone]=\zeroVec$, and $\valFunc[\state\state]=\qualFunc[\state\state]-\qualFunc[\state\ctrl]\fbPolicyHat$.
The expressions for $\ffPolicyHat$ and $\fbPolicyHat$ in~\Cref{eq:control-policy} vary depending on whether a forward or inverse dynamics formulation is used.

\subsubsection{Forward dynamics}\label{sec:forward_policy}
In forward dynamics formulations, where stagewise constraints $\consFunc(\state,\ctrl)$ are often not considered, \Cref{eq:unconstrained-policy} reduces to the classical~\gls{ddp} algorithm, i.e.,
\begin{equation}\label{eq:ddp-policy}
\begin{aligned}
\ffPolicyHat[u]\coloneqq\qualFunc[\ctrl\ctrl]^{-1}\qualFunc[\ctrl] &\quad \text{and}\quad \fbPolicyHat[u]\coloneqq\qualFunc[\ctrl\ctrl]^{-1}\qualFunc[\ctrl\state].
\end{aligned}
\end{equation}

\subsubsection{Inverse dynamics}\label{sec:inverse_policy}
In inverse dynamics formulations, when factorizing~\Cref{eq:control-policy} using the Schur complement, the feed-forward and feedback terms are computed as follows:
\begin{equation}\label{eq:schur-detailed-policy}
\begin{aligned}
\ffPolicyHat[s]\coloneqq\ffPolicy+(\ffPolicy[s]^\transpose\qualFuncEq[\ctrl\ctrl]\fbPolicyEq[s]^\transpose)^\transpose \,\,\,\text{and}\,\,\,
\fbPolicyHat[s]\coloneqq\fbPolicy+(\fbPolicy[s]^\transpose\qualFuncEq[\ctrl\ctrl]\fbPolicyEq[s]^\transpose)^\transpose
\end{aligned}
\end{equation}
with
$\ffPolicy\defeq\qualFunc[\ctrl\ctrl]^{-1}\qualFunc[\ctrl]$, $\fbPolicy\defeq\qualFunc[\ctrl\ctrl]^{-1}\qualFunc[\ctrl\state]$, $\ffPolicy[s]\defeq\consFeas-\consFunc[\ctrl]\ffPolicy$, $\fbPolicyEq[s]\defeq\qualFunc[\ctrl\ctrl]^{-1}\consFunc[\ctrl]^\transpose$, $\qualFuncEq[\ctrl\ctrl]\defeq(\consFunc[\ctrl]\qualFunc[\ctrl\ctrl]^{-1}\consFunc[\ctrl]^\transpose)^{-1}$, and $\fbPolicy[s]\defeq\consFunc[\state]-\consFunc[\ctrl]\fbPolicy$.
Alternatively, a nullspace parametrization of $\duHat$ can handle rank deficiencies in $\consFunc[\ctrl]$ both exactly and efficiently.
This approach is advantageous because it allows for parallelized computations, as detailed in \cite{mastalli-inv_dyn_mpc}.
Its resulting feed-forward and feedback terms are given by:
\begin{equation}\label{eq:nullspace-detailed-policy}
\begin{aligned}
\ffPolicyHat[n]\coloneqq\nullspaceEq\ffPolicy[n]+\qualFuncEq[\nullEqVector\nullEqVector]\fbPolicyEq[n]\consFeas \quad\text{and}\quad
\fbPolicyHat[n] = \nullspaceEq\fbPolicy[n]+\qualFuncEq[\nullEqVector\nullEqVector]\fbPolicyEq[n]\consFunc[\state]
\end{aligned}
\end{equation}
where $\nullspaceEq\in\R[\ntau\times n_Z]$ represents the nullspace basis of $\consFunc[\ctrl]$, and $\rangespaceEq\in\R[\ntau\times n_Y]$ denotes its orthonormal basis, $\ffPolicy[n]\defeq\qualFunc[\nullEqVector\nullEqVector]^{-1}\qualFunc[\nullEqVector]$, $\qualFuncEq[\nullEqVector\nullEqVector]\defeq\eyeMatrix-\nullspaceEq\qualFunc[\nullEqVector\nullEqVector]^{-1}\qualFunc[\nullEqVector\ctrl]$, $\fbPolicy[n]\defeq\qualFunc[\nullEqVector\nullEqVector]^{-1}\qualFunc[\nullEqVector\state]$, and $\fbPolicyEq[n]\defeq\rangespaceEq\consFunc[\rangeEqVector]^{-1}$.

This nullspace approach performs three decompositions for $\mathbf{Q}^{-1}_\mathbf{zz}$, $(\mathbf{h_u Y})^{-1}$, and $[\mathbf{Y} \,\, \mathbf{Z}]$.
These decompositions can be performed in parallel via efficient Cholesky and 
rank-revealing decompositions such as~\gls{lu} and QR~\cite{golub-matcompbook}.
This makes nullspace methods an attractive alternative to inexact methods. 

To compute the $\rev{\dxHat[k]}$ terms--i.e., the remaining ``primal variables'' in $\bwHat$--we perform a linear rollout as follows:
\begin{equation}\label{eq:unconstrained_linear_rollout}
\dxHat[k+1]=\dynFunc[{\state[k]}]\dxHat[k]+\dynFunc[{\ctrl[k]}]\duHat[k]+\dynFeas[k+1] \quad\text{with}\quad \dxHat[0]=\dynFeas[0],
\end{equation}
where $\duHat[k]$ defined as in~\Cref{eq:control-policy}, and this equation holds $\forall k=\{0,\cdots,N\}$.
\Cref{algo:backwarpass} summarizes the procedure for computing the endpoint-independent search direction.

\begin{algorithm}[b]
terminal value function: $\valFunc[{\state\stateN}], \valFuncHat[\stateN] = \costGrad[\state\stateN], \costGrad[\stateN]$\\
\For{$k\leftarrow N-1$ \KwTo $0$}{
compute quality function $\qualFunc$'s\hfill\Cref{eq:qual-funcs-def}\\
compute policy $\ffPolicyHat$, $\fbPolicyHat$\hfill\Cref{eq:ddp-policy,eq:schur-detailed-policy,eq:nullspace-detailed-policy}\\
compute value function $\valFuncHat[\state]$, $\valFunc[\state\state]$\hfill\Cref{eq:value-func-approx-unconstrained}
}
compute $\dxHat$, $\duHat$ via linear rollout\hfill\Cref{eq:unconstrained_linear_rollout}
\caption{Endpoint-Independent Search Direction}
\label{algo:backwarpass}
\end{algorithm}

\subsection{Endpoint-dependent search direction}
By leveraging the previously computed Riccati recursion, we can obtain the term $\bWcheck=\bA^{-1}\bB^\transpose$ efficiently.
This result stands from the following observations:
\begin{inparaenum}[(i\upshape)]
\item the Riccati recursion for $\bwHat=\bA^{-1}\ba$ partially factorizes the~\gls{kkt} system of equations,
\item this procedure exploits the Markovian structure inherent in optimal control problems, and
\item we can reuse the decomposition of the~\gls{kkt} matrix in~\Cref{eq:unconstrained-policy} to compute $\bA^{-1}\bB^\transpose$.
\end{inparaenum}
To implement this, we replace
\begin{equation}\label{eq:replace_b_by_B}
\begin{aligned}
\ba = \rev{-}
\begin{bmatrix}
\dynFeas[0] & \smallcdots & \costGrad[{\state[k]}] & \costGrad[{\ctrl[k]}] & \consFeas[k] & \dynFeas[k] & \smallcdots &  \costGrad[{\state[N]}]
\end{bmatrix}^\transpose\,\,\,\text{by}\\
\bB^\transpose = \rev{-}
\begin{bmatrix}
\zeroVec^\transpose & \smallcdots & \zeroVec^\transpose & \zeroVec^\transpose & \zeroVec^\transpose & \zeroVec^\transpose & \smallcdots & -\termFunc[{\state[N]}]^\transpose
\end{bmatrix}^\transpose.\quad
\end{aligned}
\end{equation}
This implies initializing the new Riccati recursion with $-\termFunc[{\state[N]}]^\transpose$ and solving the following system of equations:
\begin{equation}\label{eq:constrained-qual-system-inv}
\begin{aligned}
&\begin{bmatrix}
\qualFunc[\state\state] & \qualFunc[\state\ctrl] & \consFunc[\state]^\transpose \\
\qualFunc[\ctrl\state] & \qualFunc[\ctrl\ctrl] & \consFunc[\ctrl]^\transpose \\
\consFunc[\state]  & \consFunc[\ctrl] 
\end{bmatrix}
\begin{bmatrix}
\dX[\termpoint] \\ \dU[\termpoint] \\ \dL[\termpoint]
\end{bmatrix} =- 
\begin{bmatrix}
\qualFunc[\state\termpoint] \\ \qualFunc[\ctrl\termpoint] \\ \zeroVec
\end{bmatrix},
\end{aligned}
\end{equation}
where $\qualFunc[\state\termpoint]=\dynFunc[\state]^\transpose \valFuncCheck[\state\termpoint]'\in\R[\nx\times\nr]$, $\qualFunc[\ctrl\termpoint]=\dynFunc[\ctrl]^\transpose\valFuncCheck[\state\termpoint]'\in\R[\ntau\times\nr]$ as the ``Jacobians'' of the quality function associated with the endpoint-dependent search direction, and the terms $\dX[\termpoint]\in\R[\nx\times\nr]$, $\dU[\termpoint]\in\R[\ntau\times\nr]$, and $\dL[\termpoint]\in\R[\nr\times\nr]$ denote the ``primal and dual search directions'', respectively.

To exploit the problem's Markovian structure, we compute $\dU[\termpoint]$ as a function of $\dX[\termpoint]$.
This step requires solving the following condensed system:
\begin{equation}\label{eq:constrained-qual-system-inv-condensed}
\begin{aligned}
&\begin{bmatrix}
\qualFunc[\ctrl\ctrl] & \consFunc[\ctrl]^\transpose \\
\consFunc[\ctrl] 
\end{bmatrix}
\begin{bmatrix}
\dU[\termpoint] \\ \dL[\termpoint]
\end{bmatrix} =- 
\begin{bmatrix}
\qualFunc[\ctrl\termpoint] + \qualFunc[\ctrl\state]\dX[\termpoint] \\ \consFunc[\state]\dX[\termpoint]
\end{bmatrix}\scriptstyle{\in}\R[{\scriptscriptstyle{\rev{(\ntau+\nr)}\times 2\nr}}],
\end{aligned}
\end{equation}
which, for inverse dynamics, this system can be reduced to
\begin{equation}\label{eq:nullspace-constrained-kkt-reduced}
\begin{bmatrix}
\qualFunc[\nullEqVector\nullEqVector] & \qualFunc[\nullEqVector\rangeEqVector] \\
& \consFunc[\rangeEqVector]
\end{bmatrix}
\begin{bmatrix}
\dU[{\nullEqVector[\termpoint]}] \\ \dU[{\rangeEqVector[\termpoint]}]
\end{bmatrix} =-
\begin{bmatrix}
\qualFunc[{\nullEqVector\termpoint}] + \qualFunc[\nullEqVector\state]\dX[\termpoint] \\
\consFunc[\state]\dX[\termpoint]
\end{bmatrix}\scriptstyle{\in}\R[{\scriptscriptstyle{\rev{(n_Z+\nr)}\times 2\nr}}],
\end{equation}
where $\qualFunc[{\nullEqVector\nullEqVector}]\defeq \nullspaceEq^\transpose\qualFunc[\ctrl\ctrl]\nullspaceEq$, $\qualFunc[{\nullEqVector\rangeEqVector}]\defeq \nullspaceEq^\transpose\qualFunc[\ctrl\ctrl]\rangespaceEq$, $\qualFunc[{\nullEqVector\termpoint}]\defeq \nullspaceEq^\transpose\qualFunc[\ctrl\termpoint]$, $\qualFunc[{\nullEqVector\state}]\defeq \nullspaceEq^\transpose\qualFunc[\ctrl\state]$, $\consFunc[\rangeEqVector]\defeq\consFunc[\ctrl]\rangespaceEq$, and $\dU[{\nullEqVector[\termpoint]}]$, $\dU[{\rangeEqVector[\termpoint]}]$ are endpoint-dependent control policies parameterized in the null and range basis $[\mathbf{Y} \,\, \mathbf{Z}]$ obtained in~\Cref{eq:nullspace-detailed-policy}.

Both~\Cref{eq:constrained-qual-system-inv-condensed,eq:nullspace-constrained-kkt-reduced} allows \rev{us} to compute an endpoint-dependent local policy of the form
\begin{equation}\label{eq:constrained_control_policy}
\dU[\termpoint]=-\ffPolicyCheck[\termpoint]-\fbPolicyHat\dX[\termpoint],
\end{equation}
which $\fbPolicyHat$ is the same term computed in~\Cref{eq:control-policy}.

\begin{algorithm}[b]
terminal value term: $\valFunc[\state{\termpoint[N]}] = - \termFunc[{\state[N]}]^\transpose$\\
\For{$k\leftarrow N-1$ \KwTo $0$}{
compute quality terms $\qualFunc[\state\termpoint]$, $\qualFunc[\ctrl\termpoint]$\hfill\Cref{eq:constrained-qual-system-inv}\\
compute feed-forward term $\ffPolicyCheck[\termpoint]$\hfill\Cref{eq:schur-detailed-constrained-policy,eq:inverse-detailed-constrained-policy_1,eq:inverse-detailed-constrained-policy_2}\\
compute value term $\valFuncCheck[\state\termpoint]$\hfill\Cref{eq:value-func-approx-constrained}
}
compute $\dX[\termpoint]$, $\dU[\termpoint]$ via linear rollout\hfill\Cref{eq:constrained-rollout}
\caption{Endpoint-Dependent Search Direction}
\label{alg:updatebackwardpass}
\end{algorithm}

\subsubsection{Forward dynamics}
In the absence of stagewise constraints $\consFunc(\state, \ctrl)$, the feed-forward term for the endpoint-dependent control policy simplifies to:
\begin{equation}\label{eq:schur-detailed-constrained-policy}
\ffPolicyCheck[{\termpoint[u]}]\coloneqq\qualFunc[\ctrl\ctrl]^{-1}\qualFunc[\ctrl\termpoint].
\end{equation}
This simplification mirrors the result obtained in~\Cref{sec:forward_policy}.

\subsubsection{Inverse dynamics}
For inverse dynamics formulations, solving \Cref{eq:constrained-qual-system-inv-condensed} or~(\ref{eq:nullspace-constrained-kkt-reduced}) provides two methods for computing the feed-forward terms.
These methods are:
\begin{align}\label{eq:inverse-detailed-constrained-policy_1}
\ffPolicyCheck[{\termpoint[s]}]\coloneqq\ffPolicy[\termpoint]+(\ffPolicy[s\termpoint]^\transpose\qualFuncEq[\ctrl\ctrl]\fbPolicyEq[s]^\transpose)^\transpose,&\quad\quad\text{(Schur complement)}\\\label{eq:inverse-detailed-constrained-policy_2}
\rev{\ffPolicyCheck[{\termpoint[n]}]}\coloneqq\nullspaceEq\ffPolicy[n\termpoint],&\quad\quad\text{(nullspace)}
\end{align}
where $\ffPolicy[\termpoint]\defeq\qualFunc[\ctrl\ctrl]^{-1}\qualFunc[\ctrl\termpoint]$, $\ffPolicy[s\termpoint]\defeq-\consFunc[\ctrl]\ffPolicy[\termpoint]$, and $\ffPolicy[n\termpoint]\defeq\qualFunc[\nullEqVector\nullEqVector]^{-1}\qualFunc[\nullEqVector\termpoint]$.

We update the bias and ``Jacobians'' terms of the value function only as they are affected by replacing $\ba$ with $\bB^\transpose$ and we assume that $\dynFunc[\state\state]=\dynFunc[\state\ctrl]=\dynFunc[\ctrl\ctrl]=\zeroVec$.
These updates are
\begin{equation}\label{eq:value-func-approx-constrained}
\delta\valFuncCheck[\termpoint]=\frac{\Delta\valFuncCheck[{\termpoint[1]}]+2\Delta\valFuncCheck[{\termpoint[2]}]}{2}+\dX[\termpoint]^{\transpose}(\valFuncCheck[\state{\termpoint[1]}]+\valFuncCheck[\state{\termpoint[2]}])+\frac{1} {2}\dX[\termpoint]^{\transpose}\valFunc[\state\state]\dX[\termpoint]
\end{equation}
with terms computed as follows:
\begin{equation}\label{eq:value-func-cont-constrained}
\begin{aligned}
&\Delta\valFuncCheck[{\termpoint[1]}]=\ffPolicyCheck[\termpoint] ^{\transpose}\qualFunc[\ctrl\ctrl]\ffPolicyCheck[\termpoint], \quad \valFuncCheck[\state{\termpoint[1]}]=\fbPolicyHat^\transpose \qualFunc[\ctrl\ctrl]\ffPolicyCheck[\termpoint]-\qualFunc[\state\ctrl]\ffPolicyCheck[\termpoint],\\ 
&\Delta\valFuncCheck[{\termpoint[2]}]=-\ffPolicyCheck[\termpoint]^{\transpose}\qualFunc[\ctrl\termpoint], \quad\quad
\valFuncCheck[\state{\termpoint[2]}]=\qualFunc[\state\termpoint] \rev{-} \fbPolicyHat^\transpose \qualFunc[\ctrl\termpoint].
\end{aligned}
\end{equation}
For optimal control problems without stagewise constraints, these expressions simplify to $\Delta\valFunc[{\termpoint[1]}]=\ffPolicyCheck[\termpoint]^\transpose\qualFunc[\ctrl\termpoint]$, $\valFunc[\state{\termpoint[1]}]=\zeroVec$.

Finally, we compute the $\dX[\termpoint]$ terms via a linear rollout:
\begin{equation}\label{eq:constrained-rollout}
\dX[\termpoint]'=\dynFunc[\state]\dX[\termpoint] + \dynFunc[\ctrl]\dU[\termpoint], \quad\text{with}\quad \dX[{\termpoint[0]}]=\zeroVec,
\end{equation}
which provides the remaining ``primal variables'' in $\bWcheck=\bA^{-1}\bB^\transpose$.
\Cref{alg:updatebackwardpass} summarizes the procedure needed to compute the endpoint-dependent search direction.

\begin{algorithm}[t]
compute endpoint-independent direction\hfill\Cref{algo:backwarpass}\\
compute endpoint-dependent direction\hfill\Cref{alg:updatebackwardpass}\\
compute endpoint multiplier $\mulConsNext$\hfill\Cref{eq:constraint-multiplier,eq:constraint-multiplier-rengespace}\\
\ForParallel{$k\leftarrow 0$ \KwTo $N$}{
compute $\dx'$ and $\du$\hfill\Cref{eq:policy}\\
compute $\Delta\valFunc$ and $\valFunc[\state]$ \hfill\Cref{eq:value-func-cont-total}
}
\caption{Compute Direction}
\label{alg:computeDirection}
\end{algorithm}

\subsection{Endpoint multiplier and update direction}
From~\Cref{eq:saddle_multipliers}, we compute the Lagrange multiplier of the endpoint constraint $\mulConsNext$ efficiently by noting that
\begin{equation}
\begin{aligned}
\bB\bA^{-1}\bB^\transpose=\termFunc[{\state[N]}]\dX[{\termpoint[N]}],&\quad
\bB\bA^{-1}\ba=\termFunc[{\state[N]}]\dxHat[N],
\end{aligned}
\end{equation}
because $\bB$ acts as a selection matrix, $\by=\mulConsNext$, and $\bb=-\termFeas$.
Thus, the multiplier is calculated as 
\begin{equation}\label{eq:constraint-multiplier}
\mulConsNext=(\termFunc[{\state[N]}]\dX[{\termpoint[N]}])^{-1}(\termFeas+\termFunc[{\state[N]}]\dxHat[N]).
\end{equation}
Alternatively, to handle rank deficiencies in $\termFunc[{\state[N]}]$, the endpoint multiplier can be computed via 
\Cref{eq:saddle_decision_var_in_rangespace}.
This leads to
\begin{equation}\label{eq:constraint-multiplier-rengespace}
\mulConsNext=\rangespace[\termpoint](\rangespace[\termpoint]^\transpose\termFunc[{\state[N]}]\dX[{\termpoint[N]}]\rangespace[\termpoint])^{-1}\rangespace[\termpoint]^{\rev{\transpose}}(\rev{\termFeas}+\termFunc[{\state[N]}]\dxHat[N]),
\end{equation}
where $\nullspace[\termpoint]\in\R[\nr\times{\nz[c]}]$ is the nullspace basis of $\rev{\termFunc[{\state[N]}]}\dX[{\termpoint[N]}]$ and $\rangespace[\termpoint]\in\R[\nr\times{\ny[c]}]$ is chosen such that
$[\nullspace[\termpoint]\, \rangespace[\termpoint]]$ spans $\R[\nr]$.

Once $\mulConsNext$ is computed from either~\Cref{eq:constraint-multiplier} or~(\ref{eq:constraint-multiplier-rengespace}), we update the search direction $\forall k=\{0,\cdots,N\}$ as follows
\begin{equation}\label{eq:policy}
\begin{aligned}
\du=\duHat-\dU[\termpoint]\mulConsNext, &\quad\quad
\dx'=\dxHat'-\dX[\termpoint]'\mulConsNext.
\end{aligned}
\end{equation}
Moreover, the total value function is obtained by combining both Riccati passes, i.e.,
\begin{equation}\label{eq:value-func-approx-total}
\delta\valFunc=\frac{\Delta\valFunc[1]+2\Delta\valFunc[2]}{2}+\dx^{\transpose}(\valFunc[\stateone]+\valFunc[\statetwo])+\frac{1} {2}\dx^{\transpose}\valFunc[\state\state]\dx
\end{equation}
with\vspace{-2em}
\begin{align}\label{eq:value-func-cont-total}\nonumber
\Delta\valFunc[1,2]&=\Delta\valFuncHat[1,2]+\mulConsNext {^\transpose} \Delta\valFuncCheck[1,2]\mulConsNext, \quad\text{and}\\
\valFunc[{\state[1]},{\state[2]}]&=\valFuncHat[{\state[1]},{\state[2]}]+\mulConsNext {^\transpose} \valFuncCheck[{\state[1]},{\state[2]}]\mulConsNext.
\end{align}
\Cref{alg:computeDirection} summarizes the procedure to compute the total search direction.

\subsection{Expected improvement, merit function, and rollouts}
For single-shooting rollouts, as proposed in~\cite{tassa2012-synthesis}, we use~\Cref{eq:value-func-approx-total} to compute the expected improvement.
However, for feasibility-driven and multiple-shooting rollouts, we adopt a quadratic approximation of the cost function, as detailed in~\cite[Eq. (18)]{li2023-multishoot} and discussed in~\cite{mastalli2020-crocoddyl}.

These rollouts are essential for performing line search procedures.
To accept a given step length $\alpha$ during this process, we employ a merit function and update its penalty parameter $\nu$ as in~\cite{mastalli-inv_dyn_mpc}.
Our algorithm incorporates Levenberg-Marquardt regularization to enhance robustness, as in~\cite{mastalli22-boxfddp}.
For further details on the rollouts, our merit function, and their implementation, refer to~\cite{mastalli2020-crocoddyl, li2023-multishoot,martinez2024-estimation} and~\cite{mastalli-inv_dyn_mpc}.

\begin{figure}[t]
\centering
\includegraphics[width=.48\linewidth]{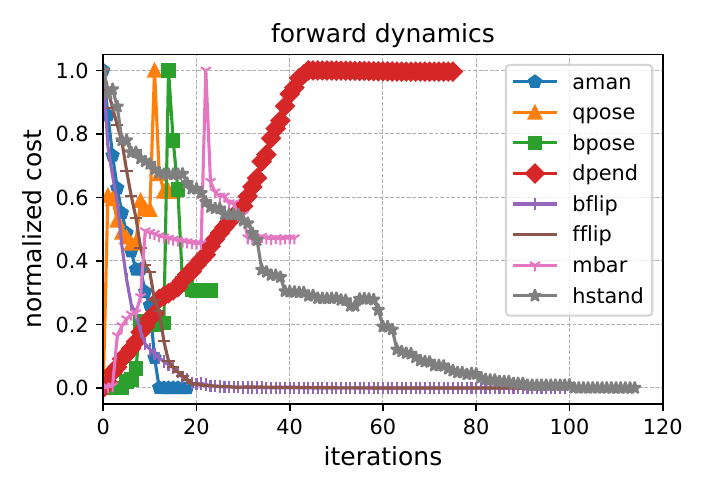}
\includegraphics[width=.48\linewidth]{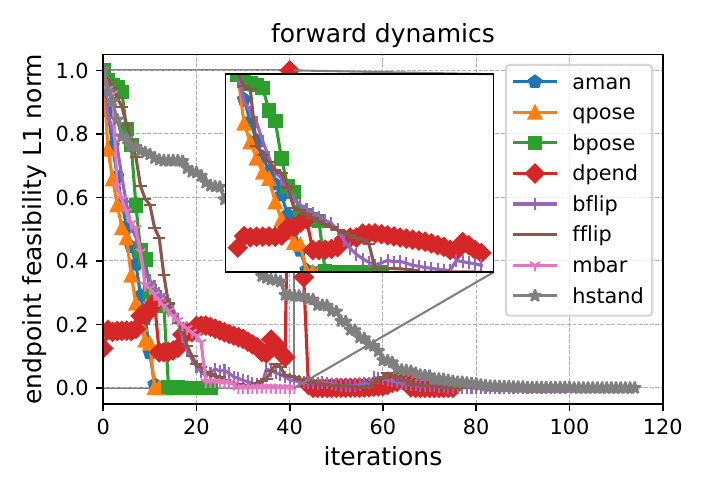}\\\vspace{-0.82em}
\includegraphics[width=.48\linewidth]{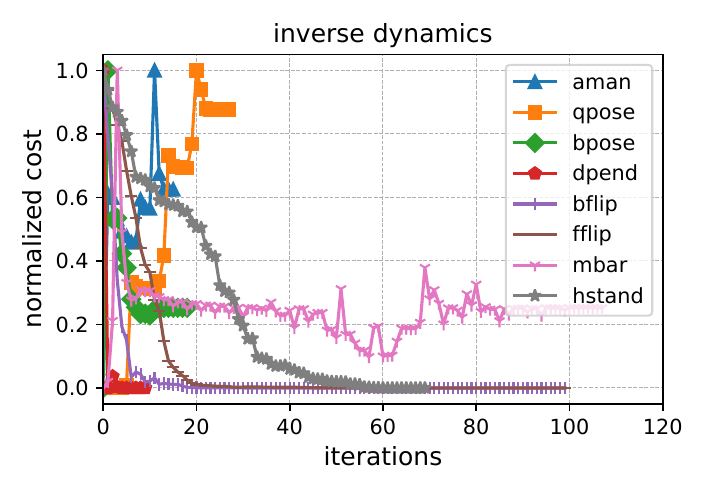}
\includegraphics[width=.48\linewidth]{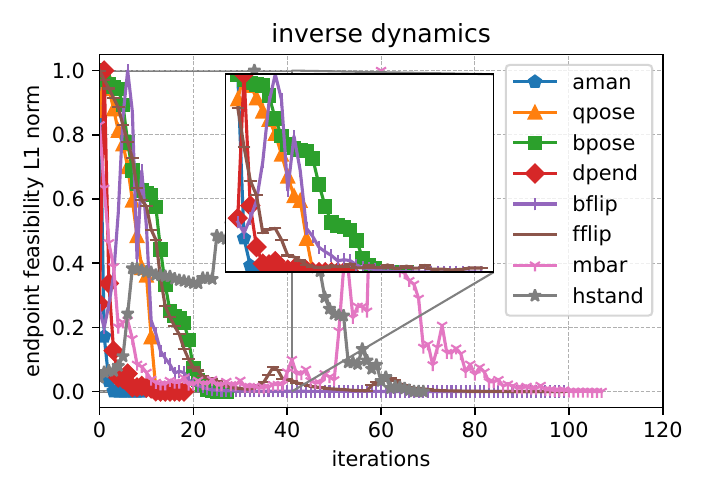}\vspace{-1em}
\caption{Cost and endpoint feasibility evolution for forward and inverse dynamics formulations:
(left) normalized cost, and (right) $\ell_1$-norm feasibility.}
\label{fig:costfeas_fwd_vs_inv}\vspace{-1em}
\end{figure}

\section{Results}\label{sec:results}
We validated our endpoint-explicit~\gls{ddp} algorithm numerically across various robotics systems (\Cref{sec:cost_feasibility_results,sec:factorization_results,sec:convergence_results}), removing running costs for path-following \rev{(keeping regularization terms only)}, as our method ensures precise endpoint satisfaction.
We evaluated:
\begin{inparaenum}[(i\upshape)]
\item aerial manipulation (\texttt{aman}), where the Borinot UAV reaches a target with its arm,
\item legged control, with the ANYmal quadruped (\texttt{qpose}) or Talos biped (\texttt{bpose}) moving toward a CoM target,
\item unstable equilibrium, where a double pendulum stabilizes at its upright position (\texttt{dpend}), and
\item gymnastics maneuvers, where the Talos humanoid performs a backflip (\texttt{bflip}), frontflip (\texttt{fflip}), handstand (\texttt{hstand}), and monkey bar (\texttt{mbar}) to specific feet or hand poses.
\end{inparaenum}
All tests were conducted on a \texttt{MacBook Pro @ Apple M2 Pro}.%

\subsection{Costs and endpoint feasibility}\label{sec:cost_feasibility_results}
To validate the effectiveness of our endpoint-explicit strategy across optimal control problems with varying stagewise constraints, we compared the cost and endpoint feasibility between forward and inverse dynamics formulations.
In~\Cref{fig:costfeas_fwd_vs_inv}, we show the normalized cost and the $\ell_1$-norm of endpoint constraint feasibility over iterations.
The results show that both cost and feasibility evolution are problem-specific.
As expected, the inverse dynamics formulation achieved faster convergence due to its ability to distribute nonlinearities more effectively throughout the optimization~\cite{ferrolho2021-inv_vs_fwd_to,mastalli-inv_dyn_mpc}.
Additionally, our nonmonotone Armijo-like condition allowed the algorithm to trade off between reducing endpoint infeasibility (e.g., \texttt{dpend}) and reducing cost (e.g., \texttt{mbar}).\vspace{-.25em}

\subsection{Endpoint factorization}\label{sec:factorization_results}
\begin{figure}[t]
\centering
\includegraphics[width=0.95\linewidth]{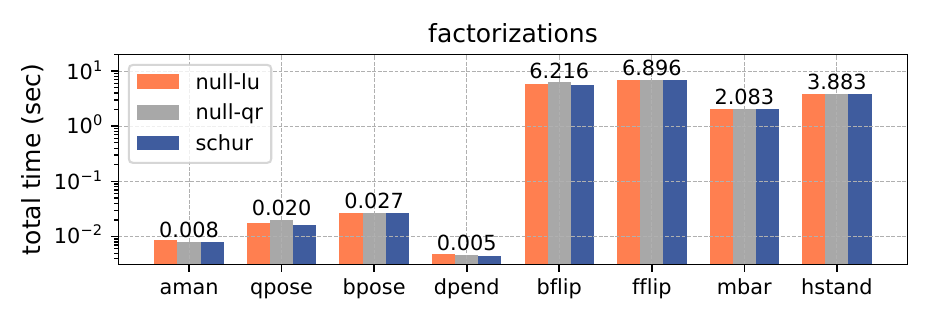}\vspace{-1em}
\caption{Average computation time for solving different factorizations and robotics problems over 10 trials, minimum time at the top.\vspace{-1.5em}}
\label{fig:computation_time}
\end{figure}

Our nullspace factorization~(\Cref{sec:abstract_nullspace}) effectively handles rank deficiencies in $\termFunc[{\state[N]}]$ while maintaining computational efficiency comparable to the traditional Schur complement method~(\Cref{sec:abstract_schur_complement}).
To evaluate this, we compared computation times between the Schur complement (\texttt{schur}) and nullspace factorizations, using~\gls{lu} with full pivoting (\texttt{null-lu}) and QR with column pivoting (\texttt{null-qr}).
\Cref{fig:computation_time} shows the average computation time for forward dynamics formulations initialized with the same warm start.
Notably, the computational complexity remained consistent across problem sizes, from $\nx+\ntau=3$ to $108$, \rev{unlike in~\cite{mastalli-inv_dyn_mpc}, where the Riccati parallelization impacted performance}. 
Therefore, the key advantage of applying our nullspace resolution is the ability to handle a broader class of problems without affecting the computational cost.\vspace{-0.25em}

\subsection{Convergence}\label{sec:convergence_results}
Evaluating convergence from random cold starts allowed us to assess our algorithm's convergence properties and demonstrate its robustness and applicability.
We conducted 10 trials with random initialization of $\stateSeq$ and $\ctrlSeq$ for both forward and inverse dynamics formulations.
The results, including the average number of iterations to achieve convergence, constraint feasibility, and success rates, are detailed in~\Cref{table:convergence-fwd-inv}.

Inverse dynamics formulations achieved a significantly higher success rate compared to forward dynamics.
For instance, in the double pendulum (\texttt{dpend}) system, the forward dynamics formulation struggled to converge and resulted in high feasibility errors.
Conversely, inverse dynamics formulations not only achieved better constraint satisfaction with less iterations but did so even as the number of constraints increased.
This is counterintuitive given its higher number of constraints compared to forward dynamics formulations.
\begin{table}[t]
\caption{Average number of iterations for convergence, constraint feasibility, and success rate over 10 trials.}
\label{table:convergence-fwd-inv}
\centering
\resizebox{0.95\linewidth}{!}{ %
\begin{tabular}{l@{\hskip 4pt}ccc@{\hskip 8pt}ccc} %
\toprule
  & \multicolumn{3}{c}{Forward Dynamics} & \multicolumn{3}{c}{Inverse Dynamics} \\
\cmidrule(lr){2-4} \cmidrule(lr){5-7}
\multirow{1}{*}{Problems} & Iter. & Feas. & Success & Iter. & Feas. & Success \\
\midrule
\texttt{qpose}  & $14.5$  & $2\mathrm{e}{-14}$   & $100\%$   & $\mathbf{15}$   & $\mathbf{2\mathrm{e}{-15}}$  & $\mathbf{100\%}$ \\
\texttt{bpose}  & $23.9$  & $6\mathrm{e}{-12}$   & $100\%$   & $\mathbf{24.3}$ & $\mathbf{2\mathrm{e}{-15}}$  & $\mathbf{100\%}$ \\
\texttt{dpend}  & $7$     & $4\mathrm{e}{+16}$   & $0\%$     & $\mathbf{46}$   & $\mathbf{1\mathrm{e}{-14}}$  & $\mathbf{100\%}$ \\
\texttt{bflip}  & $126$   & $1.3\mathrm{e}{-5}$  & $100\%$   & $\mathbf{107}$  & $\mathbf{6\mathrm{e}{-12}}$  & $\mathbf{100\%}$ \\
\texttt{fflip}  & $236$   & $3.9\mathrm{-7}$     & $100\%$   & $\mathbf{218}$  & $\mathbf{5\mathrm{e}{-10}}$  & $\mathbf{100\%}$ \\
\texttt{mbar}   & $48.7$  & $15.47$              & $0\%$     & $\mathbf{221.3}$& $\mathbf{1.5e{-3}}$  & $\mathbf{60\%}$ \\
\texttt{hstand} & $112.2$ & $15.84$              & $20\%$    & $\mathbf{69.3}$ & $\mathbf{2.70}$  & $\mathbf{50\%}$ \\
\bottomrule
\end{tabular}
}\vspace{-1.5em}
\end{table}

\subsection{Gymnastic maneuvers}\label{sec:maneuvers_results}
As discussed in~\Cref{sec:factorization_results}, our endpoint-explicit~\gls{ddp} algorithm successfully computed optimal trajectories for gymnastic maneuvers in the Talos humanoid robot.\
Despite the complexity and scale of these optimization problems, they are solved in \rev{fewer} iterations, \rev{less} seconds, and \rev{higher constraint satisfaction than penalty-based approaches such as demonstrated in~\cite{mastalli-inv_dyn_mpc}}.
\Cref{fig:snapshots} showcases Talos performing gymnastic maneuvers, including a handstand, backflip, frontflip, and monkey bar.

\subsection{MPC trials}
\begin{figure}[t]
\centering
\href{https://youtu.be/RBohdOhgbWw&t=130}{\includegraphics[width=0.4\linewidth]{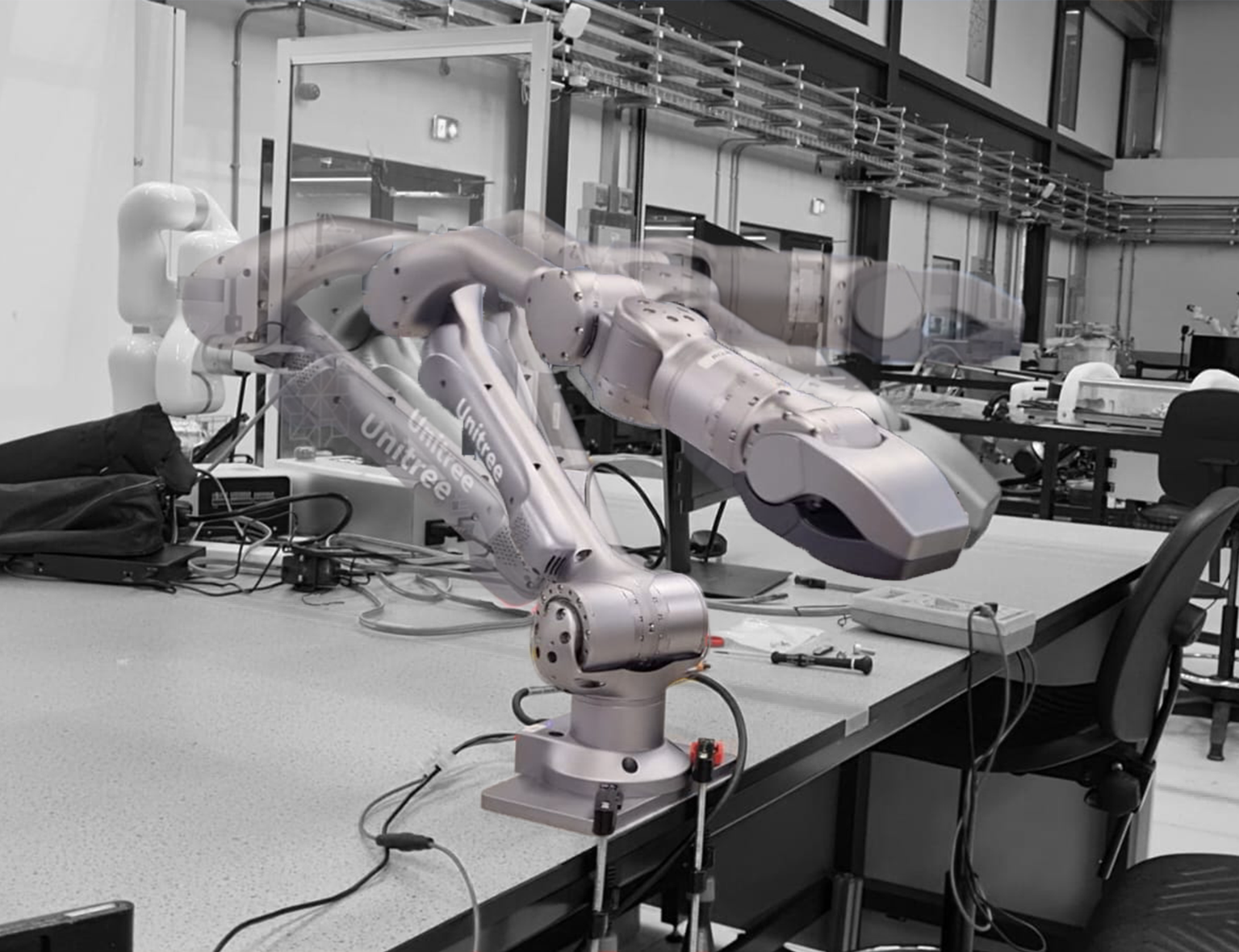}}\,\,\href{https://youtu.be/RBohdOhgbWw&t=138}{\includegraphics[width=0.55\linewidth]{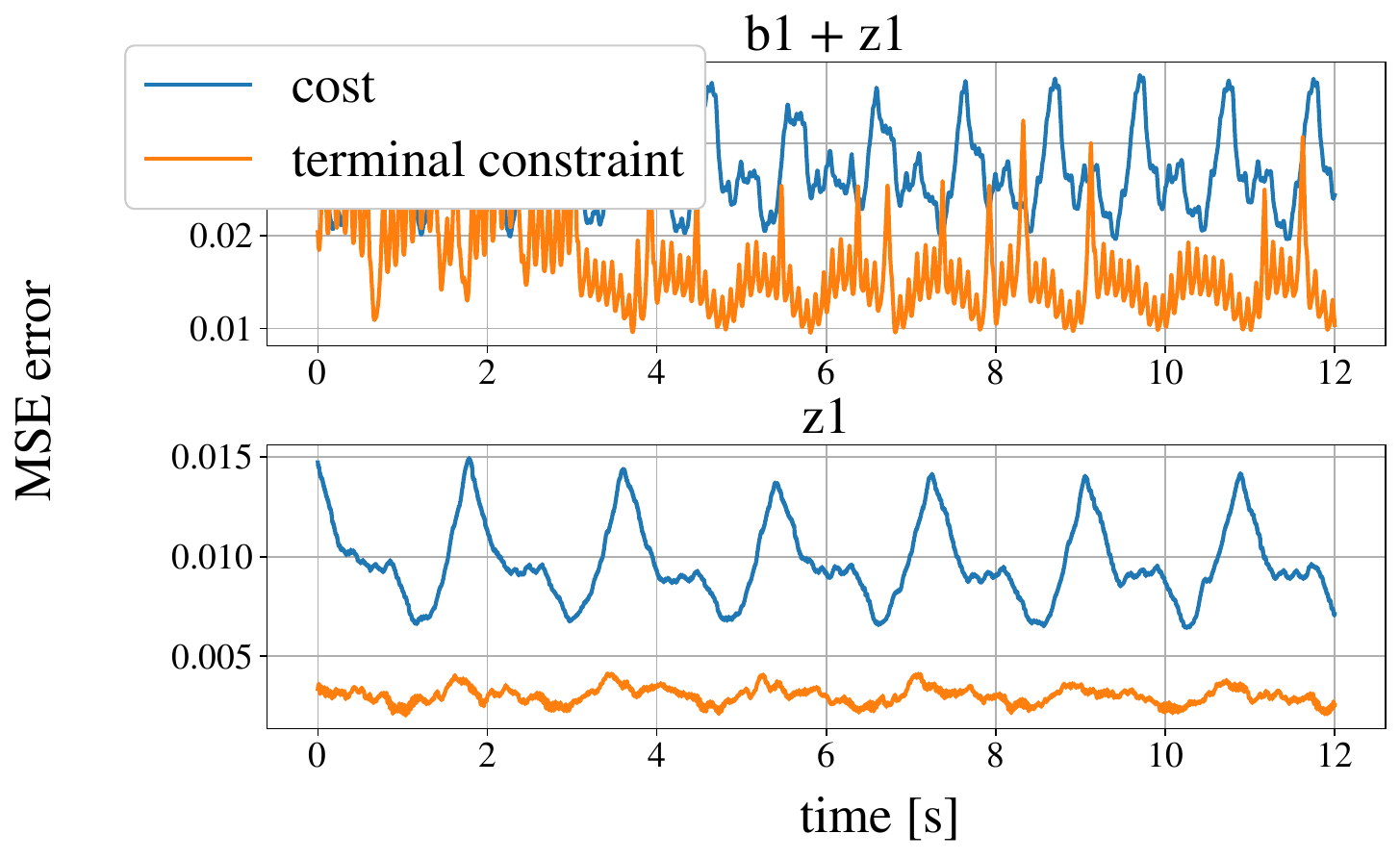}}\vspace{-0.5em}
\caption{Experiment trials with our endpoint-constrained~\gls{mpc} algorithm.
(Left) Movements of the Z1 robot demonstrating.
(Right) Tracking performance improvements in \rev{experimental trials} with the B1-Z1 quadruped robot.}
\label{fig:z1_exps}\vspace{-1.25em}
\end{figure}

To demonstrate the relevance of our endpoint-explicit method for control tasks, we conducted an experimental validation of our MPC controller.
Our method was tested in two scenarios: a simulation involving the B1 quadruped equipped with the Z1 manipulator, and a real-world experiment on the Z1 manipulator, as depicted in \Cref{fig:z1_exps} (left).
\Cref{fig:z1_exps} (right) compares the tracking error between the traditional cost-based penalty method and our approach, which incorporates a terminal constraint. In both scenarios, adding the terminal constraint led to improved tracking performance.

\section{Conclusion}
We introduced an exact method for handling endpoint constraints in~\glsdesc{ddp}, designed to efficiently tackle rank deficiencies.
Our approach demonstrated high effectiveness across a wide range of optimal control problems, adeptly handling stagewise equality constraints like inverse dynamics and contact constraints.
By reusing the Riccati recursion from problems without endpoint constraints, our method achieved exceptional computational efficiency, making it an ideal solution for real-time~\gls{mpc} in robotics.
Looking ahead, we aim to extend this strategy to manage stagewise inequality constraints and unlock its full potential for parallelizing~\gls{ddp} algorithms.

\bibliographystyle{IEEEtran}
\bibliography{references}

\end{document}